\documentclass[twoside,12pt]{article}
\usepackage{CJK}
\usepackage{indentfirst}
\usepackage{bm}
\usepackage[]{hyperref}
\usepackage[numbers,sort&compress]{natbib}
\usepackage{tikz}
\usepackage{float}
\usepackage{color}
\usepackage{graphicx}
\usepackage{subfig}
\usepackage{amsmath}
\usepackage{mathrsfs}
\usepackage{amssymb}
\usepackage{mathrsfs}
\usepackage{amsthm}
\usepackage{stfloats}

\newcommand{\RNum}[1]{\uppercase\expandafter{\romannumeral #1\relax}}

\makeatletter
\@addtoreset{equation}{section}
\makeatother
\begin{document}
\begin{CJK*}{GBK}{song}
\allowdisplaybreaks


\begin{center}
\LARGE\bf   The Global Existence  for the  Derivative Nonlinear Schr\"odinger Equation with solitons
\end{center}
\footnotetext{\hspace*{-.45cm}\footnotesize $^*$Corresponding author: Yufeng Zhang. E-mail: mathzhang@126.com, \\}
\begin{center}
\ \\Xuedong Chai, Yufeng Zhang
\end{center}
\begin{center}
\begin{small} \sl
{School of Mathematics, China University of Mining and Technology, Xuzhou, Jiangsu, 221116, People's Republic of China.
}
\end{small}
\end{center}
\vspace*{2mm}
\begin{center}
\begin{minipage}{15.5cm}
\parindent 20pt\footnotesize


 \noindent {\bfseries
Abstract}

 The  global existence of the  solution for the second-type derivative nonlinear Schr\"odinger (DNLSII) equation with solitons is presented for the first time on the line with weighted Sobolev initial data in $H^2( \mathbb{R}) \cap H^{1,1}(\mathbb{R} )$.
 The modified Darboux transformation is applied as a main research tool  to analyze the global existence for the Cauchy problem by adding or subtracting a finite number of zeros of the scattering data.
For this purpose, it is necessary to study the  invertibility and uniqueness  of the Darboux transformation first acting on the weighted Sobolev space.
 The analysis of time evolution comprehensively demonstrates the existence of a unique global solution to the DNLSII equation with Cauchy problem in $H^2( \mathbb{R}) \cap H^{1,1}(\mathbb{R} )$.

\end{minipage}
\end{center}
\begin{center}
\begin{minipage}{15.5cm}
\begin{minipage}[t]{2.3cm}{\bf Keywords}\end{minipage}
\begin{minipage}[t]{13.1cm}
 DNLSII equation; Cauchy problem; Darboux transformation; Existence of the global solution.
\end{minipage}\par\vglue8pt
\end{minipage}
\end{center}


\section{Introduction}  

In this paper, the existence of global solutions to the Cauchy problem for the  second-type derivative nonlinear Schr\"odinger (DNLSII) equation on the line is discussed.
The derivative nonlinear Schr\"odinger equation is one of the most important integrable systems in the mathematics and physics.
Chen et al. linearized the nonlinear Hamiltonian systems in order to test the integrability of the nonlinear Hamiltonian systems by the inverse scattering transformation \cite{1979DNLS}.
Several integrable equations were mentioned there, one of which is a derivative type equation
\begin{align}
& iq_t +q_{xx}+i|q|^2 q_{x}=0   ,   \label{DNLSII}
\end{align}
where the subscripts denote partial derivatives, which is the second-type derivative nonlinear Schr\"odinger equation or called Chen-Lee-Liu equation  and it is integrable.
The Global existence of DNLSII equation is studied by means of the Darboux transformation.  For this purpose, we consider the coupled DNLSII equation associated Cauchy problem
\begin{align}
& iq_t +q_{xx}+iqrq_{x}=0   ,  \quad t\in \mathbb{R}, \label{DNLSII-eq-1} \\
&  ir_t - r_{xx}+iqrr_{x}=0   ,   \label{DNLSII-eq-2}\\
& q |_{t=0}=q_0, \quad  r |_{t=0}=r_0,    \label{DNLSII-eq-3}
\end{align}
where  the initial values $q_0$ and $r_0$ fall in $H^2( \mathbb{R}) \cap H^{1,1}(\mathbb{R} )$.
Under a reduction condition $r=q^*$, the system  is reduced to the DNLSII equation.
The space $H^m( \mathbb{R}) $  represents a sobolev space with square integrability derivatives up to the order $m$. The  weighted Sobolev spaces  are defined as
\begin{align}
& L^{p,s}(\mathbb{R}) := \{ q\in L^p(  \mathbb{R}   )  :  \,\, \langle   x \rangle^s q \in L^p(\mathbb{R})      \}  , \nonumber \\
& H^{k,s}(\mathbb{R}) := \{ q\in L^{2,s}(\mathbb{R}) :  \,\,     \partial_x^j q \in L^{2,s}(\mathbb{R})   , j=1,2, \cdots, k        \} , \nonumber
\end{align}
where $\langle   x \rangle :=\sqrt{1+|x|^2}$.

With the in-depth study of the integrable systems by scholars, many methods have been applied to the  integrable equations, such as the inverse scattering method \cite{1967IST}, Hirota method \cite{1971Hirota}, B\"acklund transformation \cite{1973DT}, Darboux transformation \cite{1991DT} and so on.  Among them, the inverse scattering method is the main analytical method for the exact solution of nonlinear integrable systems and it can be used to analyze the global existence of the solution. While the $\bar{\partial}$-dressing method is closely linked to the inverse scattering method, in this aspect, we have done a lot of work on classical integrable equations \cite{2022ChaiPAMS,Chai2022AML,Chai2021TMP,Chai2023ND,Hongyi2022JMP}.

In the pioneering work \cite{1980DNLS-posedness-1,1981DNLS-posedness}, the global posedness of the Cauchy problem for the derivative nonlinear Schr\"odinger of $q_0$ in $H^2(\mathbb{R})$ is shown for initial data with a small $H^1(\mathbb{R})$ norm.
Later, Hayashi et al. \cite{1993existence,1992existence}generalized the global posedness for Cauchy condition in $H^{1}(\mathbb{R})$ with small $L^2(\mathbb{R})$ norm which is consistent with the stationary solitary waves for the DNLS equation.
The Cauchy problem of the semi classical DNLSII equation  is studied under the initial profile of rapid oscillation \cite{2002DNLS-semi-Cauchy}.

The upper bound of the $L^2(\mathbb{R})$ norm for the existence of the required initial data for the global solution has been improved by combining mass, momentum, and energy conservation with variational parameters by Wu \cite{Wu2013Global}.
A sufficient condition for global existence of solutions to a generalized derivative nonlinear Schr\"odinger equation is studied \cite{2017gDNLS-existence}.
Global well-posedness of the derivative nonlinear Schr\"odinger equation with periodic boundary condition in weighted Sobolev space $H^{\frac{1}{2}}$ is considered by Razvan \cite{2017Global-DNLS-period}.
Liu et al.\cite{2016Liu-Global-DNLS} studied the global existence of derivative nonlinear Schr\"odinger equation by direct and inverse scattering method in weighted Sobolev space. And Jenkins et al.\cite{2020Jenkins-Global-DNLS} further explored the existence of global solutions with arbitrary spectral singularities. 
The orbital stability of one-soliton solution has also been proved \cite{2006-DNLS-sta,1995Orbital-DNLS}.
The orbital stability of the sum of the two solitary waves is obtained from the variational characteristics of the medium in \cite{2017Stability-DNLS}.
Further improvements in the variational characteristics of DNLS solitaire waves and the global existence near a single solitaire wave have been explored in \cite{2017Solitary}.

In recent years, many scholars have become increasingly interested in studying the global existence problem of nonlinear evolution equations, which can be further divided into the case of with soliton and the case of without soliton.
Saalmann \cite{2017Global} combined the Riemann-Hilbert problem  with the derivative NLS equation in the presence of solitons and discussed the global existence by means of the B\"acklund transformation.
Based on the foundational and pioneering results of Deiff and Zhou \cite{2010Long},  the solvability of inverse scattering transformation is achieved.
Lipschitz continuity for direct and inverse scattering transformations has been established in appropriate spaces, indicating the global well-posedness of Cauchy problem without sharp constraints on the $L^2(\mathbb{R})$ norm of the initial data \cite{2016Existence-DNLS-IMRS}.
Simplified assumptions were made in  to exclude eigenvalues and resonances in the KN spectral problem. Excluding resonance is a natural condition for defining the so-called generic initial data with a large $L^2(\mathbb{R})$ norm \cite{2016Liu-Global-DNLS,2016Existence-DNLS-IMRS}.
If the initial datum satisfies the small norm constraint, eigenvalues are usually excluded, and in the KN spectrum problem, it is not obvious if there is an initial datum  with a large $L^2(\mathbb{R})$ norm that does not produce eigenvalues.
So far, the existence of global solutions to the nonlocal Schr\"odinger equation and  Fokas-Lenells equation on the line by Fan are studied \cite{Zhao-without,cheng2023Fokas-global}.

Furthermore, the result  is extended to the case of finite eigenvalues for spectral problems. Under the B\"acklund transform, similar to the work on the ZS spectral problem \cite{2011Deift,2014asymptotic-cNLS}, we are able to apply the inverse scattering transform technique to initial data with a finite number of solitons.
On the basis of these results, the global well-posedness with soliton  for the Cauchy problem with arbitrarily large initial data are investigated in \cite{pelinovsky2017DNLS,cheng2023Fokas-global}.

\begin{figure*}[!h]
 	\begin{center}
 	\vspace{0.8in}
 	\hspace{-0.08in}{\scalebox{0.8}[0.8]{\includegraphics{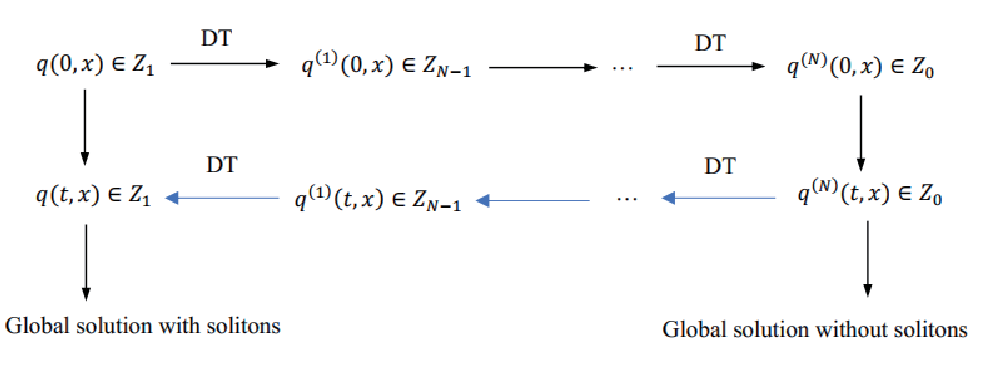}}}
 	\end{center}
 	\vspace{-0.2in} \caption{\small
 The zeros of the scattering data $a(\lambda)$ and the function $q$ can be removed or added under the Darboux transformation.}
 	\label{num1}
\end{figure*}

The goal of our present paper is to study the existence of the global solution in the case of the initial data with finite solitons in $H^{2}(\mathbb{R}) \cap H^{1,1}(\mathbb{R})$ based on the results \cite{pelinovsky2017DNLS,cheng2023Fokas-global}.
For this purpose, a modified two-fold Darboux transformation with spectral parameter $\lambda$ is constructed.
With the help of the invertibility of the Darboux transformation, we establish a bijection between the 0-soliton solution and the 1-soliton solution as shown in Fig.1.
The main results can be expressed as the following theorem.

$\mathbf{Theorem~ 1.1.}$ \quad  For every $q_0 \in H^2 (\mathbb{R}) \cap H^{1,1}(\mathbb{R})$ such that the spectral problem \eqref{Lax-DNLSII-1}   has no resonances, there is a unique global solution $q(t,\cdot) \in H^2 (\mathbb{R}) \cap H^{1,1}(\mathbb{R})$ of the Cauchy problem  for every $t\in \mathbb{R}$.

Our paper is organized as follows. In Section 2, let us briefly review some main results regarding direct  scattering method. In Section 3,  a two-fold Darboux transformation is constructed. Some properties of the Darboux transformation is presented, the invertibility of the Darboux transformation is the key to establishing 1-soliton solution and 0-soliton solution.  In Section 4, the properties of the scattering data and Jost functions under the Darboux transformation is analyzed.
In Section 5 and 6, the time evolution of the Darboux transformation and the global existence of the  solution are discussed for the DNLSII equation with soliton on the line.

\section{Construction of the scattering data and Jost functions}  

In this section, we focus on the construction of the scattering data and the Jost functions for the coupled DNLSII equation whose special case $r=\bar{q}$ is just right the scattering data and the Jost functions of the DNLSII equation. The main approach is inverse scattering method.
 The Lax pair of the coupled DNLSII equation reads \cite{1983GaugeTr}
\begin{align}
&\phi_x =  U \phi=(-i\lambda^2 +\frac{i}{4} qr)   \sigma_3 \phi +\lambda Q \phi     , \label{Lax-DNLSII-1}   \\
&\phi_t =  V \phi=(-2i\lambda^4 +i qr \lambda^2 +\frac{1}{4} (qr_x-rq_x)- \frac{i}{8}q^2 r^2 )   \sigma_3 \phi +2\lambda^3  Q \phi +\lambda W \phi ,    \label{Lax-DNLSII-2}
\end{align}
where
\begin{align}
 \sigma_3=\left(\begin{array}{cc}
1  &0 \\
0 & -1
\end{array}\right),
Q=\left(\begin{array}{cc}
0  & q \\
-r & 0
\end{array}\right),
W=\left(\begin{array}{cc}
0  & iq_x-\frac{1}{2}q^2r \\
ir_x + \frac{1}{2}r^2q & 0
\end{array}\right).
\end{align}
Taking a transformation
\begin{equation}
\psi(x,t;\lambda) = \phi(x,t;\lambda)  e^{i(\lambda^2x+2\lambda^4t) \sigma_3 } ,
\end{equation}
the Lax pair becomes
\begin{align}  \label{Lax-DNLSII-21}
&\psi_x + i\lambda^2 [\sigma_3,\psi] =( \frac{i}{4} qr\sigma_3 + \lambda Q)\psi =: \tilde{U} \psi,     \\
&\psi_t +2 i\lambda^4 [\sigma_3,\psi] = [ (iqr\lambda^2 +\frac{1}{4}(qr_x-rq_x)-\frac{i}{8}q^2 r^2)   \sigma_3  +2\lambda^3  Q  +\lambda W]\psi  =:\tilde{V} \psi,
\end{align}

Then the  Lax pair has the following asymptotic form of the Jost functions
\begin{equation}
\psi_{\pm} (x,t;\lambda)  \rightarrow I,~~x\rightarrow \pm\infty,
\end{equation}
which implies the Volterra integrable equations on the basis of integration in both directions along the parallel real axis $(-\infty,t)\rightarrow (x,t)$ and $(+\infty,t)\rightarrow (x,t)$
\begin{align}
& \psi_{-} (x,t;\lambda)  = I +\int_{-\infty}^{x} e^{-i \lambda^2 (x-y)\hat{\sigma}_{3}}\tilde{U}(y,t,\lambda) \psi_{-}(y,t,\lambda) dy,   \\
& \psi_{+}(x,t,\lambda)  = I - \int_{x}^{\infty} e^{-i \lambda^2 (x-y)\hat{\sigma}_{3}}\tilde{U}(y,t,\lambda) \psi_{+}(y,t,\lambda) dy .
\end{align}

%
The time variable $t$ is fixed and removed from the argument list of the associated function, so we omit the variable $t$ as usual. The time evolution of scattering data and the reconstruction of  the function $q$ about time $t$ are discussed later.
For two solutions of the spectral problem \eqref{Lax-DNLSII-1}, there is a transfer matrix $S(\lambda)$ based on the linear correlation of the two solutions
\begin{equation}\label{CM}
\psi_{-}(x; \lambda) e^{-i\lambda^2 \sigma_3 x} =  \psi_{+}(x; \lambda) e^{-i\lambda^2 \sigma_3 x}  S(\lambda)  .
\end{equation}
where $x\in \mathbb{R} $  is arbitrary and the scattering matrix $S(\lambda)$ is constructed by
\begin{equation}\label{S}
S(z)=\left(\begin{array}{cc}
a(\lambda)  & -\overline{b(\bar{\lambda})} \\
b(\lambda)  & \overline{ a(\bar{\lambda})}
\end{array}\right)  ,
\end{equation}
and it is easy to derive $\det S(\lambda)=1$.
By the way, the symmetry condition associated scattering matrix $S(\lambda)$ can be given as
\begin{equation}\label{sym2}
S_{\pm}(\lambda) = - \sigma \overline{S_{\pm}(\bar{\lambda})}  \sigma      .
\end{equation}
\\
$\mathbf{Proposition~ 2.1.}$
Let $q_0, r_0  \in H^2(\mathbb{R}) \cap H^{1,1}(\mathbb{R})$ and define  $\psi_{\pm}(x,t;\lambda)=(\psi_{\pm,1}(x,t;\lambda),\psi_{\pm,2}(x,t;\lambda))$ as the first and second columns of $\psi_{\pm}$.
The following properties can be found.\\
 $\bullet$  Analyticity: The Jost functions $\psi_{-,1}(x,t;\lambda)$ and $\psi_{+,2}(x,t;\lambda)$ are analytical in the first and third quadrant of the $\lambda$ plane ($D_{+}= \{  \lambda: Im\lambda^2 > 0$ \} ),  and the Jost functions $\psi_{-,2}(x,t;\lambda)$ and $\psi_{+,1}(x,t;\lambda)$ are analytical in the   second and fourth quadrant of the $\lambda$ plane ($D_{-}= \{ \lambda: Im \lambda^2 < 0  \}$).
 \\
 $\bullet$  Symmetry: The Jost functions $\psi_{\pm}(x,t;\lambda)$ admit the symmetry condition
 \begin{equation}\label{sym1}
\psi_{\pm}(x,t;\lambda) = - \sigma \overline{\psi_{\pm}(x,t;\bar{\lambda})}  \sigma      ,
\end{equation}
where
 \begin{equation}\label{s1}
\sigma=\left(\begin{array}{cc}
0  & 1 \\
-1  &  0
\end{array}\right)    .
\end{equation}
\\
 $\bullet$  Boundedness: For every $\lambda$ in $D_{+}$ and for all $q$ and $r$ satisfying $  \| \cdot\|_{L^1\cap L^{\infty}} +  \| \partial_x \cdot\|_{L^1} \leq M $ there exists a constant $C_M$ which does not depend on  $q$ and $r$, such that
 \begin{equation}\label{CM}
\| \psi_{\pm}(\cdot; \lambda)\|_{L^\infty} \leq C_M .
\end{equation}

For convenience, let us denote the vector form  $\varphi_{\pm}(x;\lambda)$    and $\phi_{\pm}(x;\lambda)$ as the first and second columns of $\psi_{\pm}$, respectively.  The analysis of the vector form retains the same properties as the matrix form, so the calculation is more concise and the calculation of the matrix results is converted into the case of the vector.
\\
$\mathbf{Proposition~ 2.2.}$ Let $q,r \in H^2(\mathbb{R}) \cap H^{1,1}(\mathbb{R})$, for every $\lambda \in \mathbb{R} \cup i \mathbb{R}$, there exist unique solutions $\varphi_{\pm}(x;\lambda) e^{-i\lambda^2x}$ and $\phi_{\pm}(x;\lambda) e^{i\lambda^2x}$ to the spectral problem  \eqref{Lax-DNLSII-1} associated with the Jost functions $\varphi_{\pm}(x;\lambda)\in L^{\infty}(\mathbb{R})$ and $\phi_{\pm}(x;\lambda)\in L^{\infty}(\mathbb{R})$ in order to satisfy the   asymptotic conditions
 \begin{equation}\label{Asycon}
 \varphi_{\pm}(x;\lambda) \rightarrow e_1, \quad     \phi_{\pm}(x;\lambda) \rightarrow e_2 ,  \quad\quad  as  \quad  x\rightarrow \pm \infty    ,
\end{equation}
as well as the scattering data, that is to say, according to the scattering relation  \eqref{CM}, the scattering coefficients $a(\lambda)$ and $b(\lambda)$ in terms  of the scattering matrix $S(\lambda)$ can be written as
\begin{align}
& a(\lambda) = \det (\varphi_{-}(x;\lambda)  e^{-i\lambda^2x}, \phi_+(x;\lambda)      e^{i\lambda^2x} ),     \label{a} \\
& b(\lambda) = \det (\varphi_{+}(x;\lambda)  e^{-i\lambda^2x}, \varphi_-(x;\lambda)     e^{-i\lambda^2x} ).  \label{b}
\end{align}

Remark: If $a(\lambda_0)=0$ for $\lambda_0 \in \mathbb{R} \cup i \mathbb{R} $, we say that $\lambda_0$ is a resonance of the spectral problem. If $a(\lambda_0)=0$ for $\lambda_0 \in \mathbb{C}_{I}:=\{   {Re(\lambda)> 0, Im(\lambda) > 0} \} $, then  $\lambda_0$ is called an eigenvalue of the spectral problem in  $\mathbb{C}_{I}$.  In particular, we call it simple if the eigenvalues satisfy the relation $a^{\prime}(\lambda_0) \neq 0$. Under the  assumption of $q,r \in L^1(\mathbb{R}) \cap L^{\infty}(\mathbb{R})$ and $q_x, r_x\in L^1(\mathbb{R})$, the scattering data $a(\lambda)$ has at most finite number of zeros in $D_+$.
Our discussion excludes resonance points but includes simple eigenvalues, and the number of eigenvalues is finite. Therefore, the initial data for the Cauchy problem may contain at most a finite number of solitons.

Denote $Z_N$ is a set of $H^2(\mathbb{R}) \cap H^{1,1} (\mathbb{R})$ such that $a(\lambda)$ has $N$ simple zeros in the field $C_{I}$, where
\begin{equation}\label{Zn}
 Z_{N} = \{q\in H^2(\mathbb{R}) \cap H^{2,1} (\mathbb{R}) ,\, a(\lambda_{j})=0,\, \lambda_j\in C_{I}, j=1,2,\ldots,N     \} .
\end{equation}
Assuming that the zero of $a(\lambda)$ is simple  for simplifying our representation. It is not a restricted assumption because of the dense in space $H^2(\mathbb{R}) \cap H^{1,1}(\mathbb{R})$.

In the following, the regularity of the Jost functions is studied.
\\
$\mathbf{Proposition~ 2.3.}$ Let $\varphi_{\pm} (x;\lambda) e^{-i\lambda^2 x}$ and $\phi_{\pm}(x;\lambda) e^{i\lambda^2 x}$ be Jost functions of the spectral problem \eqref{Lax-DNLSII-1}. For every $q,r\in  H^2(\mathbb{R}) \cap H^{2,1} (\mathbb{R}) $, $ \| \cdot\|_{H^2(\mathbb{R}) \cap H^{2,1}}\leq M $ for some $M>0$. For fixed $\lambda_{1} \in \mathbb{C}$ in $D_{+}$ , then
 \begin{align}
& \parallel \langle x \rangle \varphi_{-,2} \parallel_{L^2(\mathbb{R})} + \parallel \langle x \rangle \partial_x \varphi_{-} \parallel_{L^2(\mathbb{R})} + \parallel  \partial^2_x \varphi_{-} \parallel_{L^2(\mathbb{R})}   \leq  C_M   ,     \label{C-1} \\
&\parallel \langle x \rangle \phi_{+,1} \parallel_{L^2(\mathbb{R})} + \parallel \langle x \rangle \partial_x \phi_{+} \parallel_{L^2(\mathbb{R})} + \parallel  \partial^2_x \phi_{+} \parallel_{L^2(\mathbb{R})}   \leq  C_M      ,    \label{C-2}
\end{align}
where $\varphi_{-}:= \varphi_{-}(x;\lambda_1) = (\varphi_{-,1},\varphi_{-,2})^{T}$, $\phi_{+}:= \phi_{+}(x;\lambda_1) = (\phi_{+,1},\phi_{+,2})^{T}$ and the constant $C_M$ does not depend on $q$.
\\
\textbf{Proof}.  The statement for $\varphi_{-}$ is sufficient, and the procedure of $\phi_{+}$ is similar, so  we are just going to show the former.
Given $\varphi_- \in L^{\infty}(\mathbb{R})$ by Proposition 2.1, we want to prove $\varphi_{-,2} \in L^{2}(\mathbb{R})$, the existence of the Jost functions is uniformly in $\lambda$, it is convenient to prove the case for just $\lambda=\lambda_1$.  The integrable equation of $\varphi_{-}$ is calculated  as
\begin{equation}\label{--1}
\varphi_- = e_1+ K\varphi_-,
\end{equation}
where the operator $K$ is
\begin{equation}\label{K}
K \varphi_- = \int_{-\infty}^{x}
\left(\begin{array}{cc}
1 &0 \\
0 & e^{2i\lambda_1^2(x-y)}
\end{array}\right)
\tilde{U}(q(y))  \varphi_-(y;\lambda)  dy,
\end{equation}
which implies
\begin{equation}\label{varphi12}
\left(\begin{array}{cc}
 K \varphi_{-,1}(x;\lambda_1) \\
 K  \varphi_{-,2}(x;\lambda_1)
\end{array}\right)
=
\int_{-\infty}^{x}
\left(\begin{array}{cc}
1 &0 \\
0 & e^{2i\lambda_1^2(x-y)}
\end{array}\right)
\left(\begin{array}{cc}
  \frac{i}{4} qr & \lambda_1 q \\
 -\lambda_1 r & -\frac{i}{4}qr
\end{array}\right)
\left(\begin{array}{cc}
  \varphi_{-,1}(x;\lambda_1) \\
   \varphi_{-,2}(x;\lambda_1)
\end{array}\right) dy.
\end{equation}
For the  component of \eqref{varphi12},  it can be deduced that
\begin{align}\label{2com}
& \| \varphi_{-,2}(x,\lambda_1) \|_{L^2(-\infty,x_0)} \leq    (1+\frac{1}{2Im\lambda_1^2})  |\lambda_{1}|  \| r\|_{L^2(-\infty,x_0)}     \| \varphi_{-,1}\|_{L^\infty (-\infty,x_0)}   \nonumber\\
& \quad\quad \quad\quad \quad\quad    + \frac{1}{4}(1+\frac{1}{2Im\lambda_1^2})   \| q\|_{L^2(-\infty,x_0}   \| r\|_{L^2(-\infty,x_0)}  \| \varphi_{-,2}\|_{L^\infty(-\infty,x_0)}   .
\end{align}
Since $q$ and $r$ belong to $L^2(\mathbb{R})$, the set $\mathbb{R}$ can be divided into finitely many subintervals such that $K$ is a restriction within each sub-interval. Then we find   $\varphi_{-,2} \in L^2(\mathbb{R})$  in virtue of patching these solutions together, and we have
$  \| \varphi_{-,2}(x,\lambda_1) \|_{L^2(\mathbb{R}))} \leq C_{M}
$
where $C_{M}$ is independent of $q$ and $r$.   
From the spectral problem \eqref{Lax-DNLSII-1}, direct calculation shows that
\begin{align}\label{varphi12-1}
& \partial_x \varphi_{-,1} =   \frac{i}{4} qr  \varphi_{-,1} + \lambda_1 q \varphi_{-,2}     ,     \\
&  \partial_x \varphi_{-,2} =   - \lambda_1 r  \varphi_{-,1} -\frac{i}{4} qr  \varphi_{-,2} +2i\lambda_1^2\varphi_{-,2}    ,
\end{align}
which demonstrate that
\begin{align}\label{varphi}
& \partial_x \varphi_{-,1} \in  L^{2}(\mathbb{R})     ,    \quad     \partial_x \varphi_{-,2} \in  L^{2}(\mathbb{R})    .
\end{align}
A similar process can be obtained $\partial_x^2 \varphi_{-} \in L^2(\mathbb{R}),  x\varphi_{-}\in L^2(\mathbb{R}),  x \partial_x \varphi_{-}\in L^2(\mathbb{R}), x \partial_x^2 \varphi_{-}\in L^2(\mathbb{R})$.
$\Box$

\section{A two-fold Darboux transformation and its properties} 

In this section, a type of two-fold Darboux transformation is established in order to explore a $0$-soliton solution $q^{(1)}$ from a $1$-soliton solution $q$ in virtue of removing the zeros of scattering data $a(\lambda)$.
For this purpose, the desired Darboux transformation should be invertible. On the one hand, some modification of the Darboux transformation is made so as to satisfy the same boundary conditions for the new Jost functions. On the other hand, the existence and uniqueness for the inverse of the Darboux transformation are also discussed. 
Starting from Lax pair \eqref{Lax-DNLSII-1}-\eqref{Lax-DNLSII-2} of the coupled DNLSII equation \eqref{DNLSII-eq-1}-\eqref{DNLSII-eq-2}, a general two-fold Darboux transformation is considered by
\begin{align}\label{gDT}
\psi^{(1)} (x;\lambda)=  T \psi(x;\lambda) ,
\end{align}
where  $T=T(\lambda;\lambda_1,\lambda_2)=Q_2^{-1} T_2$, and
\begin{align}\label{gT}
&   T_2(\lambda)=\left(\begin{array}{cc}
a_2^{(2)} & 0 \\
0 & d_2^{(2)}
\end{array}\right) \lambda^2  +
\left(\begin{array}{cc}
0 & b_1^{(2)} \\
c_1^{(2)} & 0
\end{array}\right) \lambda  +
\left(\begin{array}{cc}
a_0^{(2)} & 0 \\
0 & d_0^{(2)}
\end{array}\right)    \nonumber\\
& \quad\quad\,\,   =\left(\begin{array}{cc}
a_2^{(2)}\lambda^2 + a_0^{(2)}  & b_1^{(2)} \lambda \\
 c_1^{(2)} \lambda      & d_2^{(2)}\lambda^2+d_0^{(2)}
\end{array}\right) ,   \\
&  Q_2 = [ ( \lambda_1  g_1 f_2- \lambda_2  g_2 f_1) ( \lambda_1 f_1 g_2- \lambda_2 f_2 g_1) ]^{\frac{1}{2}},
\end{align}
where each variable function can be represented as
\begin{align}\label{T-com}
& a_2^{(2)} =  \lambda_1 g_1 f_2- \lambda_2 g_2 f_1 , \quad d_2^{(2)} = \lambda_1 f_1 g_2- \lambda_2 f_2 g_1  ,  \nonumber\\
& b_1^{(2)} = -( \lambda_1^2  f_1 f_2- \lambda_2^2 f_2 f_1) , \quad c_1^{(2)} = - (\lambda_1^2  g_1 g_2- \lambda_2^2  g_2 g_1) ,  \nonumber\\
& a_0^{(2)} =  \lambda_1^2 \lambda_2  g_1 f_2- \lambda_1 \lambda_2^2  g_2 f_1 , \quad d_0^{(2)} = \lambda_1^2 \lambda_2 g_1 f_2- \lambda_1 \lambda_2^2  f_1 g_2.
\end{align}
It is worth noting  that $T_2|_{\lambda=\lambda_j} \Phi_{j}=0$ is used in the calculation of the above equation, where $ \Phi_{j}=\left(\begin{array}{cc}
 f_{j} \\
   g_{j}
\end{array}\right) $ is a vector solution of the spectral problem \eqref{Lax-DNLSII-1}-\eqref{Lax-DNLSII-2} in regard of spectral parameters $\lambda_1, \,\lambda_2 \in \mathbb{C} \backslash {\{0\}}$ and $f_{j}=f_j(x,t;\lambda), g_{j}=g_j(x,t;\lambda), j=1,2$.

Then it leads to the new  potential functions under the Darboux transformation, that is 
\begin{align}\label{q2-q}
& q^{(1)} =  \frac{\lambda_1 g_1 f_2 -\lambda_2 g_2 f_1}{\lambda_1 f_1 g_2 - \lambda_2 f_2 g_1} q -2i \frac{\lambda_1^2  f_1 f_2 -\lambda_2^2 f_2 f_1}{\lambda_1 f_1 g_2 - \lambda_2 f_2 g_1}   ,  \\
& r^{(1)} = \frac{\lambda_1 f_1 g_2 -\lambda_2 f_2 g_1}{\lambda_1 g_1 f_2 - \lambda_2 f_1 g_2} r -2i \frac{\lambda_1^2  g_1 g_2 -\lambda_2^2 g_2 g_1}{\lambda_1 g_1 f_2 - \lambda_2 f_1 g_2}   .
\end{align}

For convenience, let us denote $q^{(1)}= B(\Phi_j; \lambda) q(x,t)$  in \eqref{q2-q},  so the two-fold  Darboux transformation is given by combining \eqref{gDT} and \eqref{q2-q}.
In order to make the obtained Darboux transformation better and simpler for us to use, we need to make some modifications.

Define a bilinear form
\begin{equation}  \label{mlam}
m_{\lambda}(f,g) =  \lambda f_1 \bar{g}_1 + \bar{\lambda} f_2 \bar{g}_2   ,
\end{equation}
which acts on $\mathbb{C}^2$ with a fixed parameter $\lambda \in \mathbb{C}$.

Since the system \eqref{DNLSII-eq-1}-\eqref{DNLSII-eq-3} can be reduced to Eq.\eqref{DNLSII} when $r=\bar{q}$, we will now investigate the properties of Eq.\eqref{DNLSII}  in the following.

With the help of the above formula and the symmetry of the scattering data,
we further derive that
\begin{align}    \label{q2-DT}
 q^{(1)} =  C_{\lambda_1} (\Phi_1, \Phi_1) q -2i  D_{\lambda_1}  (\Phi_1)   =: B (\Phi_1;{\lambda_1}) q          ,
\end{align}
where
\begin{align}   \label{q2-DT-com}
C_{\lambda_1}(\Phi_1, \Phi_1) =  \frac{m_{\bar{\lambda}_1} (\Phi_1, \Phi_1) }{m_{\lambda_1}  (\Phi_1, \Phi_1)  }         ,\quad
D_{\lambda_1}(\Phi_1) = \frac{A_{\lambda_1}(\Phi_1)}{ m_{\lambda_1}(\Phi_1,\Phi_1)} ,\quad
A_{\lambda_1}(\Phi_1) = (\lambda_1^2 - \bar{\lambda}_1^2) f_1 \bar{g}_1 .
\end{align}
The Darboux transformation  \eqref{gDT} is given as
\begin{align}   \label{q2-DT-com}
\psi^{(1)}(x; \lambda) =  T (\Phi_1, \lambda, \lambda_1) \psi(x; \lambda)        ,\quad
\end{align}
associated with the Darboux matrix is
\begin{equation}\label{T}
T(\Phi_1, \lambda,\lambda_1) =\left(\begin{array}{cc}
(C_{\lambda_1}(\Phi_1) )^{\frac{1}{2}} \lambda^2 - |\lambda_1|^2(\overline{{C_{\lambda_1}(\Phi_1)}}) ^{\frac{1}{2}}     & -D_{\lambda_1}(\Phi_1) \cdot \overline{C_{\lambda_1}(\Phi_1)}^{\frac{1}{2}} \lambda \\
 -\overline{D_{\lambda_1}(\Phi_1)} \cdot (C_{\lambda_1}(\Phi_1)) ^{\frac{1}{2}} \lambda  &  (\overline{C_{\lambda_1}(\Phi_1)} )^{\frac{1}{2}} \lambda^2 - |\lambda_1|^2({C_{\lambda_1}}(\Phi_1)) ^{\frac{1}{2}}
\end{array}\right) .
\end{equation}

In the following, some modification about the Darboux transformation are made in order to satisfy the boundary conditions
\begin{equation}\label{bound}
\varphi^{(1)}(x; \lambda) \rightarrow e_1,\quad  \phi^{(1)}(x; \lambda) \rightarrow e_2 ,\quad as \quad  x \rightarrow \infty .
\end{equation}

When  $\Phi_1=\varphi_{-} e^{-i \lambda_1^2 x} $, then we have
\begin{equation}\label{bound}
T(\Phi_1, \lambda,\lambda_1)  \rightarrow T_{-\infty}(e_1, \lambda,\lambda_1)=  \left(\begin{array}{cc}
(\frac{\bar{\lambda}_1}{\lambda_1})^{\frac{1}{2}}  \lambda^2 - |\lambda_1|^2 (\frac{\lambda_1}{\bar{\lambda}_1})^{\frac{1}{2}}     & 0 \\
0    &    (\frac{\lambda_1}{\bar{\lambda}_1})^{\frac{1}{2}}  \lambda^2 - |\lambda_1|^2(\frac{\bar{\lambda}_1}{\lambda_1})^{\frac{1}{2}}
\end{array}\right) , \quad x\rightarrow -\infty,
\end{equation}
on the basis of the relation
\begin{equation}\label{D-A}
D_{\lambda_1}(\Phi_1)    \rightarrow  0,\quad    C_{\lambda_1}(\Phi_1) \rightarrow \frac{\bar{\lambda}_1}{\lambda_1}, \quad  as \quad  x \rightarrow -\infty.
\end{equation}

The new Jost functions $\varphi_{-}^{(1)}$ and $\phi_{-}^{(1)}$ are introduced as
\begin{equation}\label{newvarphi}
(\varphi_{-}^{(1)}, \phi_{-}^{(1)}) = T(\varphi_{-}, \phi_{-}) T_{-\infty}^{-1}   =(T\varphi_{-}, T\phi_{-}) \left(\begin{array}{cc}
\frac{\lambda_{1}^{\frac{1}{2}}}{ \bar{\lambda}_{1}^{\frac{1}{2}} \lambda^2-\lambda_1^2 \bar{\lambda}_1^{\frac{1}{2}}    }   & 0    \\
0  & \frac{\bar{\lambda}_{1}^{\frac{1}{2}}}{ \lambda_{1}^{\frac{1}{2}} \lambda^2-\bar{\lambda}_1^2 \lambda_1^{\frac{1}{2}}    }
\end{array}\right) .
\end{equation}

According to the asymptotics conditions \eqref{Asycon}, it is easy to find
\begin{align}    \label{varphy1}
 & \varphi_{-}^{(1)} =\frac{ \lambda_1^{\frac{1}{2}}    } {\bar{\lambda}_{1}^{\frac{1}{2}} \lambda^{2}- \lambda_1^2 \bar{\lambda}_1^{\frac{1}{2}}}  T \varphi_{-}    \rightarrow e_{1}     , \quad as \quad x \rightarrow -\infty , \\
 & \phi_{-}^{(1)} = \frac{ \bar{\lambda}_1^{\frac{1}{2}}    } {\lambda_{1}^{\frac{1}{2}} \lambda^{2}- \bar{\lambda}_1^2 \lambda_1^{\frac{1}{2}}}  T \phi_{-}    \rightarrow e_{2}     , \quad as \quad x \rightarrow -\infty .
\end{align}

Meanwhile, when  $\Phi_1=\phi_{+} e^{i \lambda_1^2 x} $ yields
\begin{equation}\label{D-A-1}
D_{\lambda_1}(\Phi_1)    \rightarrow  0,\quad    C_{\lambda_1}(\Phi_1) \rightarrow \frac{\lambda_1}{\bar{\lambda}_1}, \quad  as \quad  x \rightarrow +\infty.
\end{equation}
Naturally, we have
\begin{equation}\label{bound}
T(\Phi_1, \lambda,\lambda_1)  \rightarrow T_{+\infty}( \lambda)= \left(\begin{array}{cc}
\frac{\lambda_{1}^{\frac{1}{2}} \lambda^{2}- \bar{\lambda}_1^2 \lambda_1^{\frac{1}{2}}}{ \bar{\lambda}_1^{\frac{1}{2}}    }   & 0    \\
0  & \frac{\bar{\lambda}_{1}^{\frac{1}{2}} \lambda^{2}- \lambda_1^2 \bar{\lambda}_1^{\frac{1}{2}}}{ \lambda_1^{\frac{1}{2}}    }
\end{array}\right) , \quad x \rightarrow + \infty .
\end{equation}
Then it leads to
\begin{equation}\label{newvarphi+}
(\varphi_{+}^{(1)}, \phi_{+}^{(1)}) = T(\varphi_{+}, \phi_{+}) T_{+\infty}^{-1}   =(T\varphi_{+}, T\phi_{+}) \left(\begin{array}{cc}
\frac{ \bar{\lambda}_1^{\frac{1}{2}}    } {\lambda_{1}^{\frac{1}{2}} \lambda^{2}- \bar{\lambda}_1^2 \lambda_1^{\frac{1}{2}}}  & 0    \\
0  & \frac{ \lambda_1^{\frac{1}{2}}    } {\bar{\lambda}_{1}^{\frac{1}{2}} \lambda^{2}- \lambda_1^2 \bar{\lambda}_1^{\frac{1}{2}}}
\end{array}\right) .
\end{equation}
Furthermore, it holds the relation
\begin{align}    \label{T+}
 & \varphi_{+}^{(1)} =\frac{ \bar{\lambda}_1^{\frac{1}{2}}    } {\lambda_{1}^{\frac{1}{2}} \lambda^{2}- \bar{\lambda}_1^2 \lambda_1^{\frac{1}{2}}}  T \varphi_{+}    \rightarrow e_{1}     , \quad as \quad x \rightarrow +\infty , \\
 & \phi_{+}^{(1)} = \frac{ \lambda_1^{\frac{1}{2}}    } {\bar{\lambda}_{1}^{\frac{1}{2}} \lambda^{2}- \lambda_1^2 \bar{\lambda}_1^{\frac{1}{2}}}  T \phi_{+}    \rightarrow e_{2}     , \quad as \quad x \rightarrow +\infty .
\end{align}

Based on the above analysis, the modified Darboux matrix can be represented as
\begin{equation}\label{mT}
T(\Phi_1, \lambda,\lambda_1) = \frac{ \lambda_1^{\frac{1}{2}}    } {\bar{\lambda}_{1}^{\frac{1}{2}} \lambda^{2}- \lambda_1^2 \bar{\lambda}_1^{\frac{1}{2}}}
\left(\begin{array}{cc}
(C_{\lambda_1}(\Phi_1) )^{\frac{1}{2}} \lambda^2 - |\lambda_1|^2(\overline{{C_{\lambda_1}(\Phi_1)}}) ^{\frac{1}{2}}     & -D_{\lambda_1}(\Phi_1) \cdot \overline{C_{\lambda_1}(\Phi_1)}^{\frac{1}{2}} \lambda \\
 -\overline{D_{\lambda_1}(\Phi_1)} \cdot C_{\lambda_1}^{\frac{1}{2}} \lambda  &  (\overline{C_{\lambda_1}(\Phi_1)} )^{\frac{1}{2}} \lambda^2 - |\lambda_1|^2({C_{\lambda_1}}(\Phi_1)) ^{\frac{1}{2}}
\end{array}\right) .
\end{equation}
The variables $C_{\lambda_1}(\Phi_1)$ and $D_{\lambda_1}(\Phi_1)$ are bounded about $x$ for all $\Phi_1$, so the functions $\varphi_{\pm}^{(1)}$ and $\phi_{\pm}^{(1)}$ are also bounded in regard of $x$ for each $\lambda\in \mathbb{C} \setminus \{ \pm\lambda_{1}, \pm \bar{\lambda}_1 \} $.  Under the modified Darboux matrix, the corresponding asymptotics of the new Jost functions are
\begin{align}\label{newvarphi+-}
&  \varphi_{\pm}^{(1)} = T(e_{1}, \lambda, \lambda_{1}) e_{1} \rightarrow e_{1} , \quad as \quad  x\rightarrow \pm \infty, \\
&  \phi_{\pm}^{(1)} = T(e_{2}, \lambda, \lambda_{1}) e_{2} \rightarrow e_{2} , \quad as \quad  x\rightarrow \pm \infty.
\end{align}

Next when $x\rightarrow \pm \infty$, we will write out some  properties of $C_{\lambda}(\Phi_1)=C_{\lambda}(\varphi_{-})$  and $D_{\lambda}(\Phi_1)=D_{\lambda}(\varphi_{-})$ as follows:
\begin{align}
& C_{\lambda}(e_1) = \frac{\bar{\lambda}}{\lambda},   \,\,\,\,     C_{\lambda}(e_2) = \frac{\lambda}{\bar{\lambda}}  ,  \label{CD-pro1}  \\
&  \overline{C_{\lambda}(\Phi_1)} = C_{\bar{\lambda}} (\Phi_1), \,\,\,\,  C_{\lambda}(c \Phi_1)=C_{\lambda} (\Phi_1)        , \label{CD-pro2}\\
&  D_{\lambda} (c \Phi_1)=D_{\lambda}(\Phi_1),\,\,\,\, C_{\lambda}(\sigma_3 \Phi_1)=C_{\lambda}(\Phi_1),  \label{CD-pro3}\\
&  C_{\lambda}({\sigma_1 \Phi_1})= C_{\bar{\lambda}}(\Phi_1),\,\,\,\,  D_{\lambda}({\sigma_3 \Phi_1 })=-D_{\lambda}(\Phi_1),  \label{CD-pro4}\\
&  D_{\lambda}({\sigma_1 \Phi_1})=-\overline{D_{\lambda}(\Phi_1)},\,\,\,\,\,   D_{\lambda}(e_{1})=0=D_{\lambda}(e_2),  \label{CD-pro5}
\end{align}
where $c$ denotes a constant that is not equal to zero, these results imply the invariance of $C_{\lambda}$ and $D_{\lambda}$ by multiplying by a nonzero constant.

The following proposition  will show the potential function $q^{(1)}$   can be  obtained from different  Jost functions of the spectral problem for the coupled DNLSII equation. Also the map $B(\Phi_1, \lambda)$ is invariant as the Jost functions change.
\\
$\mathbf{Proposition~ 3.1.}    $ Let $\lambda_1 \in \mathbb{C}_{I}$  satisfy the scattering data $a(\lambda_1) =0$.   For  $q\in H^2(\mathbb{R}) \cap H^{1,1}(\mathbb{R})$, it is proved that
 \begin{align}
& q^{(1)} (x)=B(\varphi_{-}({x;\lambda_1}) e^{-i\lambda_1^2 x},{\lambda_1}) q(x)   , \label{B1} \\
&  \quad\quad \,\, \,\,\, = B(\phi_{+}({x;\lambda_1}) e^{i\lambda_1^2 x}, {\lambda_1}) q(x),    \label{B2}  \\
&  \quad\quad \,\, \,\,\, = B(\varphi_{+}({x;\bar{\lambda}_1}) e^{-i\bar{\lambda}_1^2 x}, {\bar{\lambda}_1}) q(x)  ,  \label{B3}  \\
&  \quad\quad \,\, \,\,\, = B(\phi_{-}({x;\bar{\lambda}_1}) e^{i\bar{\lambda}_1^2 x}, {\bar{\lambda}_1}) q(x)     .\label{B4}
\end{align}
\\
\textbf{Proof}.  Now let us prove each of the four equations in the above proposition one by one. The first one is obvious by means of previous discussions by setting $\Phi_1=\varphi_{-} (x;\lambda_1) e^{-i\lambda_1^2 x}$.
 As for the second equation, based on scattering data
\begin{equation}
a(\lambda_1)=\det(\varphi_{-}e^{-i\lambda_1^2 x}, \phi_{+} e^{i\lambda_1^2x}) =0 ,
\end{equation}
it intuitively obtains  $\varphi_{-}e^{-i\lambda_1^2 x}= \gamma \phi_{+} e^{i\lambda_1^2x}$, where $\gamma$ is a constant.
Then it can verify the second equation with the help of  \eqref{CD-pro1}-\eqref{CD-pro5} and $a(\lambda_1)=0$.
Now let us prove that how the fourth equation is obtained. According to the symmetry relation \eqref{sym1} and the properties of the parameter $C_{\lambda}$ and $D_{\lambda}$,  it can be proven that
\begin{equation}
C_{\lambda_1}(\varphi_{-}(x;\lambda_1)) = C_{\bar{\lambda}_1}(\sigma_1 \sigma_3 \varphi_{-}(x;\lambda_1)) = C_{\bar{\lambda}_1} (\overline{\phi_{-}(x;\bar{\lambda}_1)}) = C_{\bar{\lambda}_1} ( \phi_{-}(x;\bar{\lambda}_1)),
\end{equation}
and
\begin{equation}
D_{\lambda_1}(\varphi_{-}(x;\lambda_1)) = \overline{D_{\lambda_1}(\sigma_1 \sigma_3 \varphi_{-}(x;\lambda_1))} = \overline{D_{\lambda_1} (\overline{\phi_{-}(x;\bar{\lambda}_1)})} = D_{\bar{\lambda}_1} ( \phi_{-}(x;\bar{\lambda}_1)).
\end{equation}
Based on the map formula \eqref{q2-DT}, the forth equation is deduced.  And for the third equation, we use the conjugate version of scattering data $a(\lambda_1)=0$, that is $\varphi_{+}(x;\bar{\lambda}_1) e^{-i \bar{\lambda}_1^2 x } = \bar{\gamma} \phi_{-} (x;\bar{\lambda}_1) e^{i\bar{ \lambda}_1^2 x}$,   which leads to the third equation.   $\Box$                

The new Jost functions of the spectral problem \eqref{Lax-DNLSII-1}-\eqref{Lax-DNLSII-2} in regard of the new potential $q^{(1)}$  are expressed as
\begin{align}
&  \varphi_{-}^{(1)}(x;\lambda) =  T (\varphi_{-}(x;\lambda_1)  e^{-i \lambda_1^2 x } , \lambda, \lambda_1)  \varphi_{-}(x;\lambda)       ,   \label{nphi-1} \\
&  \varphi_{+}^{(1)}(x;\lambda) =  T (\varphi_{+}(x;\bar{\lambda}_1)  e^{-i \bar{\lambda}_1^2 x } , \lambda, \bar{\lambda}_1)  \varphi_{+}(x;\lambda)       ,\label{nphi-2}\\
&    \phi_{+}^{(1)}(x;\lambda) =  T (\phi_{+}(x;\lambda_1)  e^{i \lambda_1^2 x } , \lambda, \lambda_1)  \phi_{+}(x;\lambda)       ,\label{nphi-3}\\
&  \phi_{-}^{(1)}(x;\lambda) =  T (\phi_{-}(x;\bar{\lambda}_1)  e^{i \bar{\lambda}_1^2 x } , \lambda, \bar{\lambda}_1)  \phi_{-}(x;\lambda)       .\label{nphi-4}
\end{align}

In what follows,  the left inverse of the Darboux transformation is constructed, and the one-to-one correspondence are established between $1$-soliton and $0$-soliton solution.
\\
$\mathbf{Proposition~ 3.2.} $    The left inverse of the map $B(\Phi_1; \lambda_1)$ exists uniquely.
\\
\textbf{Proof}.  Assuming the left inverse exists and is denoted as $B(\tilde{\Phi}_1; \lambda_1)$, then the following formula holds
\begin{align}    \label{qtilde}
&  \tilde{q}= B (\tilde{\Phi}_1; \lambda_1)  B(\Phi_1; \lambda_1)  q       ,     \nonumber \\
&  \,\,\,\,  = C_{\lambda_1}(\tilde{\Phi}_1) C_{\lambda_1}(\Phi_1) q -2i C_{\lambda_1}(\tilde{\Phi}_1) D_{\lambda_1} (\Phi_1)      -2i D_{\lambda_1} (\tilde{\Phi}_1),
\end{align}
where
\begin{align}
& C_{\lambda_1}(\tilde{\Phi}_1) C_{\lambda_1}(\Phi_1)=1        ,   \label{q-3} \\
& C_{\lambda_1}(\tilde{\Phi}_1) D_{\lambda_1} (\Phi_1)  +  D_{\lambda_1} (\tilde{\Phi}_1) =0   . \label{q-4}
\end{align}
So we can calculate that
\begin{align}
& C_{\lambda_1}(\Phi_1) =\overline{C_{\lambda_1}(\tilde{\Phi}_1) }        ,   \label{q-3} \\
& D_{\lambda_1}(\tilde{\Phi}_1) =   -\overline{C_{\lambda_1}(\Phi_1) } D_{\lambda_1} (\Phi_1)    . \label{q-4}
\end{align}
On account of  two expressions mentioned above and the fact that the relation $C_{\lambda_1}(\Phi_1)=\overline{C_{\lambda_1} (\tilde{\Phi}_1)}$ can be equivalently written as the following equation
\begin{align}
|f_1|^2 |\tilde{f}_1|^2=|g_1|^2|\tilde{g}_1|^2        ,   \label{fgmodel}
\end{align}
which implies that
\begin{align}
|\tilde{f}_1|=n |g_1|,\quad\quad      |\tilde{g}_1|=n |f_1|    ,   \label{fgmodel-1}
\end{align}
where $n$ is a positive number that satisfies Eq.\eqref{fgmodel}.
From the perspective of the plural, it is evident that
\begin{align}
\tilde{f}_1=n_1  \bar{g}_1,\quad\quad      \tilde{g}_1=n_2 \bar{f}_1    ,   \label{fgmodel-2}
\end{align}
with $n_1,\,n_2\in \mathbb{C}$ are constants and  $|n_1|=|n_2|$.

Consequently, the left inverse of $B(\Phi_1, \lambda_1)$ given by \eqref{fgmodel-2} exists associated with $\tilde{\Phi}_1$.


Next  the uniqueness of the variable $\tilde{\Phi}_1$ can be represented  by function $\bar{\Phi}_1$, that is to say, let $\tilde{\Phi}_1=(\tilde{f}_1,\tilde{g}_1)^T$ be given by
\begin{align}   \label{spe1}
\tilde{f}_{1}    = \frac{\bar{g}_1}{m_{\lambda_1}(\Phi_{1},\Phi_{1})} - \frac{\bar{g}_1}{m_{\bar{\lambda}_1}(\Phi_{1},\Phi_{1})},  \quad
\tilde{g}_{1}    = \frac{\bar{f}_1}{m_{\bar{\lambda}_1}(\Phi_{1},\Phi_{1})} - \frac{\bar{f}_1}{m_{\lambda_1}(\Phi_{1},\Phi_{1})}.
\end{align}
In fact, based on the spectral problem \eqref{Lax-DNLSII-1}, there is
\begin{align}   \label{spe1-Lax}
\partial_x \Phi_{1}    =  (-i \lambda_1^2  + \frac{i}{4}  |q|^2) \sigma_3 \Phi_{1}  +\lambda_1 Q \Phi_{1}  ,
\end{align}
According to the expression \eqref{mlam} and above equation, we can easily get
\begin{align}  \label{lam1xx}
& \partial_x m_{\lambda_1}(\Phi_1,\Phi_1)    =  (\lambda_1^2 - \bar{\lambda}_1^2)( -i\lambda_1 |f_1|^2 + i \bar{\lambda}_1 |g_1|^2 +q \bar{f}_1 g_1)            ,        \nonumber \\
& \partial_x m_{\bar{\lambda}_1}(\Phi_1,\Phi_1)    =  (\bar{\lambda}_1^2 - \lambda_1^2) ( i \bar{\lambda}_1 |f_1|^2 - i \lambda_1 |g_1|^2 +\bar{q}  f_1 \bar{g}_1),
\end{align}
the  formula can be further received
\begin{align}  \label{spec1}
& \tilde{f}_{1,x} + i\lambda_1^2 \tilde{f}_1 -\frac{i}{4} |q^{(1)}|^2 \tilde{f}_1
=  \frac{\lambda_1 \bar{f}_1}{ m_{\lambda_1}^2} [-\bar{\lambda}_1 q |f_1|^2 - \lambda_1 q |g_1|^2 + 2i(\lambda_1^2-\bar{\lambda}_1^2) \bar{g}_1 f_1]     \nonumber \\
& -\frac{\lambda_1 \bar{f}_1}{ m_{\lambda_1} m_{\bar{\lambda}_1}^2}  [- \bar{\lambda}_1^2 |f_1|^4 q - 2 \bar{\lambda}_1 \lambda_1 |f_1|^2 |g_1|^2 q - \lambda_1^2 |g_1|^4q  + 2i \bar{\lambda}_1(\lambda_1^2-\bar{\lambda}_1^2) |f_1|^2 f_1\bar{g}_1 +2i \lambda_1 (\lambda_1^2-\lambda_1^2) f_1\bar{g}_1 |g_1|^2]
\nonumber \\
& = \frac{\lambda_1 \bar{f}_1}{ m_{\lambda_1}^2}  [- q m_{\bar{\lambda}_1} + 2i(\lambda_1^2 - \bar{\lambda}_1^2)   \bar{g}_1 f_1]
+\frac{\lambda_1 \bar{f}_1}{ m_{\lambda_1} m_{\bar{\lambda}_1}^2}  [-m_{\bar{\lambda}_1}^2 q + 2im_{\bar{\lambda}_1} (\lambda_1^2-\bar{\lambda}_1^2)  f_1 \bar{g}_1]   \nonumber \\
& =\lambda_1 q^{(1)} \tilde{g}_1,
\end{align}
it is noted  that the constraint is
\begin{align}  \label{restr}
\frac{|\lambda_1|^2 |f_1|^4 (\lambda_1^2-\bar{\lambda}_1^2)}{m_{\bar{\lambda}_1}^2} = \frac{1}{4} (|q|^2-|q^{(1)}|^2).
\end{align}
And the $\tilde{g}_1$ in the spectral problem \eqref{spe1-Lax} can be similarly deduced such that \eqref{spe1} is proven.
$\Box$

\section{Application of the Darboux transformation}   

In this section, the smoothness of the new Jost functions is discussed based on the modified Darboux transformation. And the zeros of the scattering data are subtract so as to present that the new potential $q^{(1)}\in Z_0$ when the old potential  $q \in Z_1 $. On the basis of these results, the regularity of the new potentials under the Darboux transformation is proven.
\\
$\mathbf{Proposition~ 4.1.}$  Given that $\varphi_{\pm}^{(1)}$ and $\phi_{\pm}^{(1)}$ are defined by \eqref{nphi-1}-\eqref{nphi-4},  then their removeable singularities are $\lambda=\pm\lambda_1$ and $\lambda=\pm\bar{\lambda}_1$ in their analytic domian, respectively.
\\
\textbf{Proof}.
Considering that one of the Jost functions $\varphi_{-}^{(1)}(x;\lambda)$ is sufficient and the other case follow similarly. Let $\varphi_{-}=(\varphi_{-,1}, \varphi_{-,2})^{T}$ and $\varphi_{-}^{(1)}=(\varphi_{-,1}^{(1)}, \varphi_{-,2}^{(1)})^{T}$ be
2-vectors, we can calculated for $\lambda \in C_{I} \cup C_{III} \setminus \{ \pm \lambda_1\}$,
\begin{eqnarray}   \label{varphi-11}
&&  \varphi_{-,1}^{(1)}(\lambda)  = \frac{\lambda_1^{\frac{1}{2}}}{ \bar{\lambda}_1^{\frac{1}{2}} \lambda^2 - \lambda_1^2 \bar{\lambda}_1^{\frac{1}{2}}}\cdot  \left[ \left(\frac{m_{\bar{\lambda}_1}(\varphi_{-},\varphi_{-})}{ m_{\lambda_1}(\varphi_{-},\varphi_{-})}\right)^{\frac{1}{2}}   \lambda^2  \varphi_{-,1}(\lambda) - |\lambda_1|^2 (\frac{m_{\lambda_1}(\varphi_{-},\varphi_{-})}{m_{\bar{\lambda}_1}(\varphi_{-},\varphi_{-})})^{\frac{1}{2}}  \varphi_{-,1}(\lambda)-\right.   \nonumber\\
&&  \left.\frac{ (\lambda_1^2- \bar{\lambda}_1^2) \varphi_{-,1}(\lambda_1) \overline{\varphi_{-,2}(\lambda_1)}}{ m_{\lambda_1}(\varphi_{-}, \varphi_{-})} (\frac{m_{\lambda_1}(\varphi_{-},\varphi_{-})}{m_{\bar{\lambda}_1}(\varphi_{-},\varphi_{-})})^{\frac{1}{2}} \cdot  \lambda  \cdot \varphi_{-,2}(\lambda)      \right]   .
\end{eqnarray}
In order to obtain the desired result, the following transformation is performed as
\begin{align}  \label{varphi-12}
&  \left(\frac{m_{\bar{\lambda}_1}(\varphi_{-},\varphi_{-})}{ m_{\lambda_1}(\varphi_{-},\varphi_{-})}\right)^{\frac{1}{2}} \varphi_{-,1}^{(1)}(\lambda)    = \left(\frac{\lambda_1}{ \bar{\lambda}_1}\right)^{\frac{1}{2}}  \frac{1}{ \lambda^2 - \lambda_1^2}\cdot   \left[ \lambda^2\cdot \frac{m_{\bar{\lambda}_1}(\varphi_{-},\varphi_{-})}{ m_{\lambda_1}(\varphi_{-},\varphi_{-})} \varphi_{-,1}(\lambda) - |\lambda_1|^2 \varphi_{-,1}(\lambda) - \right.  \nonumber\\
&  \left. \frac{\lambda^2 - \bar{\lambda}_1^2 }{m_{\lambda_1}(\varphi_{-},\varphi_{-}) \lambda} \varphi_{-,2} (\lambda)      \right] = \left(\frac{\lambda_1}{ \bar{\lambda}_1}\right)^{\frac{1}{2}}  \frac{(\lambda^2 - \lambda_1^2) \cdot \bar{\lambda}_1  |\varphi_{-,1}(\lambda_1)|^2\varphi_{-,1}(\lambda)+H(\lambda)    }{ (\lambda^2 - \lambda_1^2) m_{\lambda_1}(\varphi_{-},\varphi_{-})},
\end{align}
where
\begin{align}
H (\lambda) = (\lambda^2-\bar{\lambda}_1^2) \lambda_1 | \varphi_{-,2}(\lambda_1)|^2  \varphi_{-,1}(\lambda) -\lambda (\lambda_1^2-\bar{\lambda}_1^2) \varphi_{-,1}(\lambda_1) \overline{\varphi_{-,2}(\lambda_1)}  \varphi_{-,2}(\lambda).
\end{align}
In the light of odevity of the above equation, in other words, $\varphi_{-,1}$ is even about $\lambda$ and $\varphi_{-,2}$ is odd about $\lambda$, so it is obvious that $H(x; \lambda_1)= H(x;-\lambda_1)=0$.
According to Proposition 2.1,  $H(\lambda)$ is analytic in $C_{I}\cup C_{III}$, hence, there exists $H(\lambda)= (\lambda^2-\lambda_1^2) \tilde{H}(\lambda)$, where $\tilde{H}(\lambda)$ is also analytic in $C_{I}\cup C_{III}$. Further calculations can determine
\begin{align}  \label{varphi-12}
& \varphi_{-,1}^{(1)}(\lambda)    =  \left(\frac{\lambda_1}{ \bar{\lambda}_1}\right)^{\frac{1}{2}}    \frac{  \bar{\lambda}_1 |\varphi_{-,1}({\lambda}_1)|^2 \varphi_{-,1}(\lambda) +\tilde{G}(\lambda) } {m_{\lambda_1}(\varphi_{-},\varphi_{-})}         \left(\frac{m_{\lambda_1}(\varphi_{-},\varphi_{-})}{ m_{\bar{\lambda}_1}(\varphi_{-},\varphi_{-})}\right)^{\frac{1}{2}}           ,
\end{align}
which explicates that $\pm\lambda_1$ are removable singularities of $\varphi_{-,1}^{(1)}(\lambda)$. Similar process can be used to prove that
 $\varphi_{-,2}^{(1)}(\lambda)$ has removable singularities $\pm\lambda_1$.   $\Box$
 \\
$\mathbf{Proposition~ 4.2.}$
The relation between $\tilde{\Phi}_1$ and the new Jost function is expressed as
\begin{align}  \label{Phi-11}
& \tilde{\Phi}_1 (x)    = (\frac{1}{\bar{\lambda}_1}-\frac{1}{\lambda_1})         \frac{1 }{ \gamma a^{(1)}(\lambda_1) }  e^{-i\lambda_1^2 x} \varphi_{-}^{(1)} (x;\lambda_1) + (\frac{1}{\bar{\lambda}_1}-\frac{1}{\lambda_1}) \frac{1 }{ a^{(1)}(\lambda_1) }  e^{i\lambda_1^2 x} \phi_{+}^{(1)} (x;\lambda_1),
\end{align}
where $\lambda_1$ is fixed such that $a(\lambda_1)=0$ as well as $a^{\prime}(\lambda_1)\neq 0$, $\gamma \neq 0$  is a constant governed by $a(\lambda_1)=0$ and $a^{\prime}(\lambda_1)\neq 0$.
\\
\textbf{Proof}.
Let 2-vectors $\tilde{\Phi}_1$  and $\varphi_{-}$ be of the form $\tilde{\Phi}_{1}=(\tilde{f}_{1}, \tilde{g}_{1})^{T}$ and $\varphi_{-}=(\varphi_{-,1}, \varphi_{-,2})^{T}$, respectively.  Starting from the components of $\tilde{\Phi}_{1}$, that is to say,
\begin{align}
& \tilde{f}_1 (x)    = \frac{e^{i\lambda_1^2x} \overline{\varphi_{-,2}(x;\lambda_1)}}{m_{\lambda_1}(\varphi_{-}(x;\lambda_1), \varphi_{-}(x;\lambda_1) )} -  \frac{e^{i\lambda_1^2x} \overline{\varphi_{-,2}(x;\lambda_1)}}{m_{\bar{\lambda}_1}(\varphi_{-}(x;\lambda_1), \varphi_{-}(x;\lambda_1) )}
 ,  \label{Phi-12} \\
&   \tilde{g}_1 (x)    = \frac{e^{i\lambda_1^2x} \overline{\varphi_{-,1}(x;\lambda_1)}}{m_{\bar{\lambda}_1}(\varphi_{-}(x;\lambda_1), \varphi_{-}(x;\lambda_1) )} -  \frac{e^{i\lambda_1^2x} \overline{\varphi_{-,1}(x;\lambda_1)}}{m_{\lambda_1}(\varphi_{-}(x;\lambda_1), \varphi_{-}(x;\lambda_1) )} , \label{Phi-13}
\end{align}
because of the boundary condition
\begin{align}
\lim_{x\rightarrow -\infty} \varphi_{-}(x;\lambda_1)   = e_1   , \nonumber
\end{align}
we establish
\begin{align} \label{lim1}
\lim_{x\rightarrow -\infty} e^{- i \lambda_1^2 x } \tilde{\Phi}_1 = (\frac{1}{\bar{\lambda}_1}- \frac{1}{\lambda_1})  e_2   .
\end{align}
When the case of $a(\lambda_1)=0$, there is a constant $\gamma$ such that
\begin{align} \label{r}
\varphi_{-}(x;\lambda_1)  e^{-i\lambda_1^2x}    =   \gamma  \phi_{+}(x;\lambda_1)  e^{i\lambda_1^2x}  ,  \quad x\in \mathbb{R}.
\end{align}
Combining  with the  formula \eqref{Phi-12}-\eqref{Phi-13} yields
\begin{align}
& \tilde{f}_1 (x)    =  \frac{e^{-i\lambda_1^2x} \overline{\phi_{+,2}(x;\lambda_1)}}{\gamma m_{\lambda_1}(\phi_{+}(x;\lambda_1), \phi_{+}(x;\lambda_1) )} -  \frac{e^{-i\lambda_1^2x} \overline{\phi_{+,2}(x;\lambda_1)}}{ \gamma m_{\bar{\lambda}_1}(\phi_{-}(x;\lambda_1), \phi_{-}(x;\lambda_1) )}
,  \label{Phi-14} \\
&   \tilde{g}_1 (x)    =\frac{e^{-i\lambda_1^2x} \overline{\phi_{+,1}(x;\lambda_1)}}{\gamma m_{\bar{\lambda}_1}(\phi_{+}(x;\lambda_1), \phi_{+}(x;\lambda_1) )} -  \frac{e^{-i\lambda_1^2x} \overline{\phi_{+,1}(x;\lambda_1)}}{\gamma m_{\lambda_1}(\phi_{+}(x;\lambda_1), \phi_{+}(x;\lambda_1) )}    . \label{Phi-15}
\end{align}
According to the other boundary condition
\begin{align}
\lim_{x\rightarrow +\infty} \phi_{+}(x;\lambda_1)   = e_2   , \nonumber
\end{align}
it can be seen that
 \begin{align} \label{lim2}
\lim_{x\rightarrow +\infty} e^{i \lambda_1^2x} \tilde{\Phi}_1 =\frac{1}{\gamma}(\frac{1}{\bar{\lambda}_1} - \frac{1}{\lambda_1}) e_1       .
\end{align}

It is known that $\tilde{\Phi}_1$ is a solution to the spectral problem \eqref{Lax-DNLSII-1} with new potentials $q^{(1)}$ for the parameter $\lambda=\lambda_1$, based on the analyticity, any solution can be represented as
 \begin{align}  \label{Phi-111}
& \tilde{\Phi}_1 (x)    = \varsigma_1   e^{-i\lambda_1^2 x} \varphi_{-}^{(1)} (x;\lambda_1)  + \varsigma_2  e^{i \lambda_1^2x} \phi_{+}^{(1)}(x;\lambda_1)        ,
\end{align}
where $\varsigma_1 $ and $\varsigma_2 $ remains to be determined and is independent  of $x$. Further analysis and calculation are manipulated
 \begin{align}  \label{Phi-112}
& \lim_{x\rightarrow -\infty} e^{-i \lambda_1^2x} \tilde{\Phi}_1 = \varsigma_2  a^{(1)}(\lambda_1) e_2 , \\
& \lim_{x\rightarrow +\infty} e^{i \lambda_1^2x} \tilde{\Phi}_1 = \varsigma_1  a^{(1)}(\lambda_1) e_1       .
\end{align}
Note that when $a^{(1)} (\lambda_1) \neq 0$, $\varsigma_1$  and $\varsigma_2$ can be carried out
 \begin{align}  \label{xishu}
\varsigma_{1} = \frac{1}{\gamma}(\frac{1}{\bar{\lambda}_1} - \frac{1}{\lambda_1}) \frac{1}{a^{(1)}(\lambda_1)}
 ,\quad\quad    \varsigma_{2} =      (\frac{1}{\bar{\lambda}_1} - \frac{1}{\lambda_1}) \frac{1}{a^{(1)}(\lambda_1)}                 ,
\end{align}
Therefore, the above process demonstrates the proposition 4.2.   $\Box$

Remark:  \eqref{CD-pro1}-\eqref{CD-pro5} have already described the invariance of Darboux transformation even though $\Phi_1$ is multiplied by a nonzero constant. Consequently, the formula in Proposition 4.2 can be rewritten as
 \begin{align}  \label{4.3-1}
& \tilde{\Phi}_1 (x)    =   e^{-i\lambda_1^2 x} \varphi_{-}^{(1)} (x;\lambda_1) + \alpha_1 e^{ i\lambda_1^2 x}    \phi_{+}^{(1)}(x;\lambda_1)        .
\end{align}

On the basis of  the definition of $Z_N$ given by \eqref{Zn},  it satisfies $a(\lambda_1)=0$ when $\lambda_1 \in C_{I} $ is the root of $a(\lambda)$, so there is $q \in Z_1 $.  The scattering data $a(\lambda)$ and $b(\lambda)$ are transformed so as to present that the new potential $q^{(1)}\in Z_0$ when the old potential  $q \in Z_1 $.
\\
$\mathbf{Proposition~ 4.3.}$  There exists  $\lambda_1 \in C_{I} $ and $q \in Z_1 $ satisfy $a(\lambda_1)=0$ when $\Phi_1(x)=\varphi_{-} (x;\lambda_1) e^{-i\lambda_1^2 x}$, where $\varphi_{-}$ is the Jost function of the spectral problem \eqref{Lax-DNLSII-1}. Then, the potential under the Darboux transformation satisfy $q^{(1)} = B(\Phi_1,\lambda_1) q \in Z_0$ .
\\
\textbf{Proof}.  Assuming that $a(\lambda)$ has a simple zero $\lambda=\lambda_1$, under the Darboux transformation, the  new scattering datum becomes
\begin{align}
a^{(1)}(\lambda) = \det(\varphi_{-}^{(1)}(x, \lambda), \phi_{+}^{(1)}(x, \lambda))    ,   \label{a1=0}
\end{align}
has no zeros in $C_{I}$, where Jost functions $\varphi_{-}^{(1)}$ and $\phi_{+}^{(1)}$ has been shown by \eqref{nphi-1}-\eqref{nphi-4}. Further  calculations are
\begin{align}  \label{a1=0-1}
& a^{(1)}(\lambda)   = \det(T [\varphi_{-} (x, \lambda_1)  e^{-i\lambda_1^2 x }, \lambda,\lambda_1] \varphi_{-}(x; \lambda), \,  T[ \phi_{+} (x, \lambda_1) e^{i\lambda_1^2 x}, \lambda, \lambda_1] \phi_{+}  (x; \lambda)  ),   \nonumber \\
& \quad\quad \quad   = \det(T [\varphi_{-} (x, \lambda_1)  , \lambda,\lambda_1] \varphi_{-}, \,    T[ \phi_{+} (x, \lambda_1) , \lambda, \lambda_1] \phi_{+}),      \nonumber \\
& \quad\quad \quad   = \det(T [\varphi_{-} (x, \lambda_1)  , \lambda,\lambda_1] \varphi_{-}, \,    T[ \varphi_{-} (x, \lambda_1)  , \lambda, \lambda_1] \phi_{+}),      \nonumber \\
& \quad\quad \quad   =  a(\lambda)  \det(T [\varphi_{-} (x, \lambda_1)  , \lambda,\lambda_1] ) ,    \nonumber \\
& \quad\quad \quad   =  a(\lambda)  \frac{\lambda_1}{\bar{\lambda}_1} \frac{\lambda^2 - \bar{\lambda}_1^2}{\lambda^2 - \lambda_1^2}    .
\end{align}
It can be analyzed that $a(\lambda)$ has no zeros point on $C_{I}$ because $\lambda_1$ is a only simple zero point of $a(\lambda)$.

On the flip side, the transformation from $b(\lambda)$ to $b^{(1)}$ are calculated similarly as
\begin{align}  \label{b1=0-1}
& b^{(1)}(\lambda)   = \det(  e^{-2i\lambda^2 x } \varphi_{+}^{(1)}(x; \lambda) , \,  \varphi_{-}^{(1)}(x; \lambda) ) ,   \nonumber \\
& \quad\quad\quad  = \det(  e^{-2i\lambda^2 x } \varphi_{+}^{(1)}(x; \lambda) , \,      T[ \varphi_{-} (x;\lambda_1) e^{-i\lambda_1^2 x} , \lambda, \lambda_1]      \varphi_{-}(x; \lambda) ) ,   \nonumber \\
& \quad\quad\quad  = \det(  e^{-2i\lambda^2 x } \varphi_{+}^{(1)}(x; \lambda) , \,      T[ \phi_{+} (x;\lambda_1) e^{i\lambda_1^2 x} , \lambda, \lambda_1]      [a(\lambda) \varphi_{+}(x; \lambda)   + e^{2i\lambda^2 x}b(\lambda) \phi_{+}(x;\lambda)  ] ) ,   \nonumber \\
& \quad\quad\quad  =   b(\lambda)\det(e_1, T(e_2, \lambda, \lambda_1)e_2)  ,  \nonumber \\
& \quad\quad\quad   =b(\lambda)  .
\end{align}
It can be seen that $b(\lambda)$ has not added new zeros and singularities.  $\Box$

Next the transformation of the new potentials are constructed.
\\
$\mathbf{Proposition~ 4.4.}$  
Under the same conditions as in Proposition~ 4.3, for each $q^{(1)} \in H^{2}(\mathbb{R})  \cap H^{1,1}(\mathbb{R}) $, there is $ \| q^{(1)}\|_{  H^{2}(\mathbb{R})  \cap H^{1,1}(\mathbb{R}) } \leq M$ for some constant $M > 0$, the transformation
 \begin{align}  \label{4.4}
 B(\tilde{\Phi}_1 ,\lambda_1) q^{(1)} \in   H^{2}(\mathbb{R})  \cap H^{1,1}(\mathbb{R})       ,
\end{align}
satisfies
 \begin{align}  \label{4.4-1}
 \| B(\tilde{\Phi}_1 ,\lambda_1) q^{(1)} \|_{  H^{2}(\mathbb{R})  \cap H^{1,1}(\mathbb{R})  }   \leq C_{M}  ,
\end{align}
where $C_M$ is a constant and independent of $q^{(1)}$.
\\
\textbf{Proof}.
Reviewing the expression (\ref{a}), we find that
 \begin{align}  \label{a1lambda}
 & |a^{(1)}(\lambda_1)| =  \mid  (\varphi_{-,1}^{(1)} e^{-i \lambda_1^2 x} +\alpha_1 \phi_{+,1}^{(1)} e^{ i\lambda_1^2 x}) \phi_{+,2}^{(1)} e^{i\lambda_1^2 x}  \nonumber\\
 &  \quad\quad\quad\quad\quad      -( \varphi_{-,2}^{(1)} e^{-i\lambda_1^2 x} + \alpha_1 \phi_{+,2}^{(1)} e^{i\lambda_1^2 x}  ) \phi_{+,1}^{(1)} e^{i\lambda_1^2 x} \mid     \nonumber \\
 &    \quad\quad\quad\quad\quad   \leq  \|  \phi_{+}^{(1)} (\cdot;\lambda_1) \|_{L^{\infty}} \cdot ( |e^{i\lambda_1^2 x }\tilde{f}_1|+|e^{i\lambda_1^2 x} \tilde{g}_1|    )       ,
\end{align}
owing to $a^{(1)}(\lambda_1) \neq 0$ and $|m_{\lambda}(\Phi_1,\Phi_1) | \geq | Re(\lambda_1)| (|f_1|^2 + |g_1|^2)$, there exists a constant $C_{M}> 0$ that does not depend on the new potential $q^{(1)}$, which leads to
 \begin{align}  \label{Cm}
 &  \frac{1}{    |m_{\lambda_1}(e^{i\lambda_1^2 }\tilde{\Phi}_1(x), e^{i\lambda_1^2 }\tilde{\Phi}_1(x)) |}
 \leq     C_{M}  .
\end{align}
Analogous discussion can be done, that is to say
 \begin{align}  \label{a1lambda-1}
 &  | a^{(1)}(\lambda_1) |  \leq  |\alpha_1|^{-1}  \| \varphi_{-}^{(1)}(\cdot; \lambda_1)\|_{L^{\infty}}\|   ( |e^{-i\lambda_1^2 x }\tilde{f}_1|+|e^{-i\lambda_1^2 x } \tilde{g}_1|    )       ,
\end{align}
such that $\frac{1}{    |m_{\lambda_1}(  e^{-i\lambda_1^2 x }    \tilde{\Phi}_1(x), e^{- i\lambda_1^2 x } \tilde{\Phi}_1(x)) |}
 \leq     C_{M}$.

For the bound,  we have
 \begin{align}  \label{D1}
 &  \left|\frac{\tilde{f}_1 \tilde{\bar{g}}_1}{m_{\lambda_1}(\tilde{\Phi}_1,\tilde{\Phi}_1)}\right|  \leq   \frac{ \left|  \varphi_{-,1}^{(1)} (x;\lambda_1)   \overline{\varphi_{-,2}^{(1)} (x;\lambda_1)}   \right| } { \left| m_{\lambda_1}(e^{-i\lambda_1^2 x }\Phi_1, e^{-i\lambda_1^2 x }\Phi_1 )  \right| }
  +   |\alpha_1|^2 \frac{\left| \phi_{+,1}^{(1)}(x;\lambda_1)  \overline{\phi_{+,2}^{(1)}(x;\lambda_1)} \right| }     { \left| m_{\lambda_1} (e^{i\lambda_1^2 x }    \tilde{\Phi}_1(x),e^{i\lambda_1^2 x }    \tilde{\Phi}_1(x))  \right| }   \nonumber\\
  &  +  |\alpha_1| \frac{ \left|  \varphi_{-,1}^{(1)} (x;\lambda_1)  \overline{\phi_{+,2}^{(1)}(x;\lambda_1)} \right| +  \left| \phi_{+,1}(x;\lambda_1)^{(1)}     \overline{\varphi_{-,2}^{(1)}(x;\lambda_1) }    \right|   }{  \left| m_{\lambda_1}(\tilde{\Phi}_1,\tilde{\Phi}_1)   \right| }  ,
\end{align}
recalling the bounds \eqref{C-1} and \eqref{C-2}, we obtain
 \begin{align}  \label{D2}
 &  \| D_{\lambda_1} (\Phi_1) q^{(1)}\|_{L^{2,1}} \leq C_M      ,
\end{align}
based on the similar proof of  Proposition~3.2, it can be seen that
\begin{align}  \label{D3}
 &  \| B(\Phi_1, {\lambda_1}) q^{(1)} \|_{L^{2,1}} \leq C_M      .
\end{align}
Therefore, $  \|  \partial_x (B(\Phi_1, {\lambda_1}) q^{(1)}) \|_{L^{2,1}} $ and $  \|  \partial_x^2 (B(\Phi_1, {\lambda_1}) q^{(1)}) \|_{L^{2}} $ can also be estimated similarly from Proposition~2.3.  So the bound  \eqref{4.4-1} of the Proposition 4.4 is completed. $\Box$
\\
$\mathbf{Lemma~ 4.5.}$    For a initial value problem of a system of linear ordinary differential equations with the initial value problem
\begin{align}
 & \frac{dy}{dx}= A(x) y    ,  \label{ODE1} \\
 & y(x_0)=0 ,\quad\quad   x_0 \in [\alpha, \, \beta],   \label{ODE2}
\end{align}
where $y=(y_{1}(x), \cdots,y_{n}(x) )^{T}$ and  $A(x)$ is an $n$-order matrix-value function. On account of the standard Picard approximation, the above system associated with the initial value problem exists a unique zero solution.

On the basis of the regularity of $\varphi(x;\lambda)$ and $\phi(x;\lambda)$ proven in Proposition 2.3, we discover that the Darboux transformation can be given by an operator from $q\in H^{2} (\mathbb{R})\cap H^{1,1}(\mathbb{R}) $ to $q^{(1)} \in H^{2} (\mathbb{R})\cap H^{1,1}(\mathbb{R}) $. It can be expressed as the following proposition.
\\
$\mathbf{Proposition~ 4.6.}$
 For the potential  $q\in H^{2} (\mathbb{R})\cap H^{1,1}(\mathbb{R}) $, there exists
 $q^{(1)} \in H^{2} (\mathbb{R})\cap H^{1,1}(\mathbb{R}) $, where we have defined $\Phi_1(x)=\varphi_{-} e^{-i \lambda_1^2 x}$ and $\varphi_{-}$ is the Jost solution of the spectral problem \eqref{Lax-DNLSII-1}.
\\
\textbf{Proof}.
we found that
\begin{align}
q^{(1)}=\pm q,
\end{align}
by calculating $C_{\lambda}= \pm 1$ and $D_{\lambda}=0$ when the case of $\lambda \in \mathbb{R} \cup i \mathbb{R}$. Therefore, the transformation  has no sense, so we concentrate on when $\lambda$ is outside the continuous spectrum.
With the help of the Lemma 4.5 and notice that $Re (\lambda_1) > 0$, we have
\begin{align}
m_{\lambda_1}(\Phi_1,\Phi_1)=0  \Leftrightarrow \Phi_1=0,
\end{align}
thus, $\Phi_1^{\prime} (x_0)=0$ according to $\Phi_1 (x_0)=0$ at $x_0 \in \mathbb{R}$, which derives $\Phi_1(x)=0$  for every $x\in \mathbb{R}$. It is obvious to observe that $\Phi_1(x) = \varphi_{-} (x;\lambda_1) e^{-i\lambda_1^2 x} \neq 0$ and $m_{\lambda_1}(\Phi_1,\Phi_1) \neq 0$ for every finite $x\in \mathbb{R}$. For convenience, we use $m_{\lambda_1} (\varphi_-,\varphi_-)$   instead of $m_{\lambda_1}(\Phi_1,\Phi_1) $.

When $a(\lambda_1) \neq 0$, there is a constant $\gamma_1$ such that
\begin{align}  \label{4.6-1}
|m_{\lambda_1}(\varphi_-, \varphi_-)| \geq \gamma_1,
\end{align}
for all $x\in \mathbb{R}$. That is to say,
\begin{align}
\lim_{x\rightarrow -\infty} m_{\lambda_1}(\varphi_-, \varphi_-)= \lambda_1.
\end{align}
According to the bound $ \phi_{+}(x, \lambda_1)  \rightarrow e_2$ as $x\rightarrow +\infty$  and $\varphi_{-}(x; \lambda_1) \in L^\infty (\mathbb{R})$, we have
\begin{align}
 \lim_{x\rightarrow +\infty} \varphi_{-}(x;\lambda_1) =a(\lambda_1) \neq 0,
\end{align}
which implies
\begin{align}
\lim_{x\rightarrow +\infty} m_{\lambda_1}(\varphi_-, \varphi_-) \neq 0.
\end{align}
Hence, the inequality \eqref{4.6-1}  holds. Making use of the trigonometric inequalities and $|C_{\lambda_1} (\Phi_1)| =1$, the following equation can be estimated
\begin{align}   \label{q1norm}
& \| q^{(1)}\|_{L^{2,1}} \leq \| q\|_{L^{2,1}} +2\| D_{\lambda_1}(\Phi_1)\|_{L^{2,1}} \nonumber\\
&  \quad\quad\quad\quad    \leq   \| q\|_{L^{2,1}} + 4 \gamma_1 ^{-1} |\lambda_1^2-\bar{\lambda}_1^2| \| \varphi_{-,1}(\cdot;\lambda_1) \overline{\varphi_{-,2} (\cdot; \lambda_1)}\|_{L^{2,1}}  < \infty    .
\end{align}
As for the norms $\| \partial_x q^{(1)}\|_{L^{2,1}}$ and $\| \partial_x^2 q^{(1)}\|_{L^{2}}$ can be estimated similarly.

On the other hand, when $a(\lambda_1) =0$, Eq.\eqref{4.6-1} is invalid because of
\begin{align}
\lim_{x\rightarrow +\infty} m_{\lambda_1}(\varphi_-, \varphi_-) = 0.
\end{align}
The norm of $q^{(1)}$ in $L^{2,1}$ are estimated on the interval $(-\infty, R)$, where $R> 0$ is arbitrary. In order to extend the estimate \eqref{q1norm} to the whole interval $(R, +\infty)$, taking advantage of the relation
\begin{align}
\varphi_-(x;\lambda_j)  e^{-i\lambda_j^2 x} = \gamma \phi_+(x;\lambda_j)  e^{i\lambda_j^2 x}, \quad x\in \mathbb{R}.
\end{align}
Combining with the \eqref{B2}, there is
\begin{align}
\lim_{x\rightarrow +\infty} m_{\lambda_1}(\phi_+, \phi_+) = \bar{\lambda}_1,
\end{align}
with the help of the equivalent representation of $q^{(1)}$, the estimate on the interval $(R,+\infty)$ can also be done similarly.  $\Box$

\section{Time evolution of the Darboux transformation}  

It is known that $\varphi_{-}(t,x;\lambda)$ is time-dependent according to linear spectral problem \eqref{Lax-DNLSII-1}-\eqref{Lax-DNLSII-2}.
Now let us elaborate on the idea of $1$-soliton solution and extend the argument to the case of finitely many solitons iteratively. Let $q(t, \cdot) \in Z_1 \subset H^{2}(\mathbb{R}) \cap H^{1,1}(\mathbb{R})$ be a local solution of the initial problem on $(-T, T)$ for some $T>0$.
For the  fixed time $t\in (-T,T)$, there is a new potential $q^{(1)}(t,\cdot) = B(\Phi_1(t,\cdot), \lambda_1) q(t,\cdot)$ under the Darboux transformation of the spectral problem \eqref{Lax-DNLSII-1}.
When $\lambda_1$ belongs to $C_{I}$, we have $q^{(1)}(t,\cdot) \in Z_0$ if $a(\lambda_1)=0$.
In addition, let $\tilde{q} \in Z_0$ be a solution to the spectral problem with the Cauchy problem and the initial value is $\tilde{q} (0, \cdot) = q^{(1)}(0,\cdot) \in Z_0$, which implies that $\tilde{q} (t, \cdot) \in Z_0$ for $t\in(-T,T)$.
The Fig.2 can illustrate the idea.
\begin{figure*}[!h]
 	\begin{center}
 	\vspace{0.8in}
 	\hspace{-0.08in}{\scalebox{0.8}[0.8]{\includegraphics{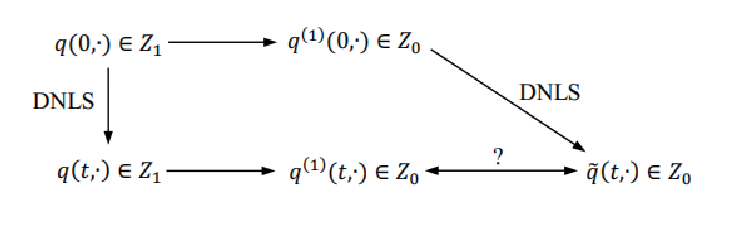}}}
 	\end{center}
 	\vspace{-0.2in} \caption{\small
 Diagram of the relationship between different solutions.}
 	\label{num2}
\end{figure*}
Therefore, the key is  to prove that $\tilde{q}(t,\cdot)=q^{(1)}(t,\cdot)$ for every $t \in (-T,T)$. For this purpose, we first prove that the two potentials generate the same scattering data.
\\
$\mathbf{Lemma~ 5.1.}$
For every $t \in (-T,T)$, the potentials $\tilde{q}(t,\cdot)$ and $q(t,\cdot)$ have the same scattering data .
\\
\textbf{Proof}.
The both potentials $\tilde{q}(t,\cdot)$ and  $q(t,\cdot)$ are in $Z_0$ for every $t \in (-T,T)$.
For the potential $q(t,\cdot) \in Z_1$ with $t \in (-T,T)$, we have $ l(x;\lambda) = \frac{b(t,\lambda)}{a(t,\lambda)}$ for $\lambda \in \mathbb{R}\cup i\mathbb{R}$. Then, denote $l^{(1)}(t,\lambda) = \frac{b^{(1)}(t,\lambda)}{a^{(1)}(t,\lambda)}$,  there is the reflection coefficients of $q^{(1)}(t,\cdot) \in Z_0$ for $t \in (-T,T)$.
On the basis of the Proposition 4.3, the relation of the new and old reflection coefficient  can be expressed as
\begin{align}  \label{lab}
l^{(1)}(t, \lambda) =  l(t,\lambda) \frac{\bar{\lambda}_1}{\lambda_1}  \frac{\lambda^2- \lambda_1^2}{\lambda^2-\bar{\lambda}_1^2}    ,    \quad    \lambda \in \mathbb{R}\cup i\mathbb{R} , \quad t\in(-T,T).
\end{align}
When $q(x,t)$ is a solution to DNLSII equation in $Z_1$, the time evolution of the reflection coefficient $l(t,\lambda)$ is derived as
\begin{align}  \label{lab1}
l(t, \lambda) =  l(0,\lambda)  e^{-4i \lambda^4 t} , \quad t\in(-T,T).
\end{align}
At the same time, the reflection coefficient $l^{(1)}(t,\lambda)$ is
\begin{align}  \label{lab2}
l^{(1)}(t, \lambda) =  l^{(1)}(0,\lambda) e^{-4i \lambda^4 t}   , \quad t\in(-T,T).
\end{align}
The reflection coefficient $\tilde{l}$ of the potential $\tilde{q}$, there is $l^{(1)}(0,\lambda) = \tilde{l}(0,\lambda)$ because of $q^{(1)}(0,\cdot) = \tilde{q}(0,\cdot)$. Therefore, the reflection coefficient is represented as
\begin{align}  \label{lab3}
\tilde{l} (t, \lambda) = \tilde{l} (0; \lambda) e^{-4i \lambda^4 t}  =  l^{(1)}(0,\lambda) e^{-4i \lambda^4 t} = l^{(1)}(t,\lambda)   , \quad t\in(-T,T).
\end{align}
So the process demonstrates the Lemma.   $\Box$
\\
$\mathbf{Proposition~ 5.2.}$
The potential $q^{(1)}(t,\lambda) = B(\Phi_1(t,\lambda), \lambda_1) q(t,\lambda)$ is a new solution for the DNLSII equation in $t \in (-T,T)$ .
\\
\textbf{Proof}.
The existence and Lipschitz continuity of the map
$
L^{2,1} (\mathbb{R}\cup i\mathbb{R}) \supset X \ni l \longmapsto q\in Z_0 \subset H^2(\mathbb{R}) \cap H^{1,1}(\mathbb{R})
$
can be established with the help of the associated Riemann-Hilbert problem. The map is bijective as well as  $\tilde{q}(t,\lambda)=q^{(1)}(t,\lambda)$ for $t\in (-T,T)$ under Lemma 5.1. Therefore, $q^{(1)}$ is a solution of DNLSII equation based on $\tilde{q}$ is a solution.  $\Box$

\section{Global existence with solitons} 

Based on the Lemma 5.1 and Proposition 5.2, we conclude that
\\
$\mathbf{Proposition~ 6.1.}$
For fixed $\lambda_1 \in C_{I}$, on the basis of a local solution $q(t, \lambda) \in H^2(\mathbb{R}) \cap H^{1,1}(\mathbb{R})$ for the Cauchy problem for some $T> 0$, let
\begin{align}  \label{Phi-6.1}
\Phi_1 (t,x) := \varphi_{-}(t,x,\lambda_1) e^{-i(\lambda_1^2 x+ 2\lambda^4 t)}
,
\end{align}
where $\varphi_{-}(t,x,\lambda_1)$ is also the Jost solution of the linear system \eqref{Lax-DNLSII-1}-\eqref{Lax-DNLSII-2}. So there is $q^{(1)}(t,\cdot) = B(\Phi_1(t,\cdot), \lambda_1) q(t,\cdot) \in H^{2}(\mathbb{R}) \cap H^{1,1}(\mathbb{R})$ for $t\in [0,T)$ and $q^{(1)}$ satisfies the Cauchy problem for $q^{(1)} (0,\cdot)= B(\Phi_1(0,\cdot)) q(0,\cdot)$.
\\
\textbf{Proof}.
The proof in the case of finitely many solitons depends on the iteration of the above discussion. The distinct eigenvalues $\{ \lambda_1, \lambda_1,\cdots,\lambda_N\}$ in $C_{I}$ are removed by iterating the Darboux transformation $N$ times for a given $q\in Z_N,\, N\in \mathbb{N}$. Define $q^{(0)} =q$ and
\begin{align}  \label{Phi-6.1-1}
q^{(s)}=B(\Phi_1^{(s-1)}, \lambda_1 ) q^{(s-1)}  , \quad   1\leq s \leq N ,
\end{align}
which implies that $q^{(N)} \in Z_0$. The proof of Lemma 5.1 and Proposition 5.2 apply to the last potential $q^{(N)}$. Consequently, the $N$-fold iteration of the Darboux transformation in regard of  a solution $q(t,\cdot) \in Z_N$  with $y\in[0,T)$ generate a new solution $q^{(N)} (t,\cdot) \in Z_0$ for the initial problem.
 The iteration diagram is shown in the Fig.1.  $\Box$
\\
\textbf{Proof of the Theorem 1.1}.
Let $q_0 \in Z_1\subset H^{2}(\mathbb{R}) \cap H^{1,1}(\mathbb{R})$ and $\lambda_1 \in C_{IV}$ be the only root of $a(\lambda)$ in $C_{I}$.
 According to the Proposition~ 4.3, there is $q_0^{(1)} = B(\Phi_1,\lambda_1)q_0 \in Z_0 \subset H^{2}(\mathbb{R}) \cap H^{1,1}(\mathbb{R})$ if $\Phi_1(x)=\varphi_- (x;\lambda_1) e^{-i\lambda_1^2 x}$, where $\varphi_- $ is the Jost solution of the spectral problem about $q_0$.
Recalling Proposition~3.2, the map is invertible with $q_0 = B(\tilde{\Phi}_1, \lambda_1) q_0^{(1)}$, where $\tilde{\Phi}_1$ has been expressed by means of the new Jost functions $\varphi_{-}^{(1)}$ and $\phi_{+}^{(1)}$  from the Proposition 4.2.
let $T > 0$ be the maximal existence time for the solution $q(t,\cdot) \in Z_1,\, t\in (-T,T)$ in terms of the initial value $q_0 \in Z_1$ and the eigenvalue $\lambda_1$.
The solution $q(t,\cdot) \in Z_1$ has the Jost functions $\{ \varphi_{\pm}(t,x,\lambda), \phi_{\pm}(t,x,\lambda) \}$ for fixed $t\in (-T,T)$. In this case, $q^{(1)}$ is defined by the B\"acklund transformation
\begin{align}  \label{Phi-6.1-2}
q^{(1)}=B(\Phi_1(t,\cdot), \lambda_1 ) q  , \quad   \Phi_1(t, \cdot)= \varphi_{-} (t,x,\lambda_1) e^{-i\lambda_1^2 x- 2i\lambda_1^4 t} ,
\end{align}
where the definition of $\varphi_{-} (t,x,\lambda_1) $ is similar to the previous discussion mentioned in this article.

The solution $q^{(1)}(t,\cdot) \in Z_0$ for $ t\in(-T,T)$ is a solution of the Cauchy problem with the initial value $q_0^{(1)} \in Z_0$ according to Proposition 5.2.
It can be analyzed that the solution $q^{(1)}\in Z_0$ is uniquely continued for $t\in \mathbb{R}$.  Similar to the solution $q(t,x)$, $q^{(1)}$ admits the Jost functions $\{ \varphi_{\pm}^{(1)}(t,x,\lambda), \phi_{\pm}^{(1)}(t,x,\lambda) \}$.  For $t\in (-T,T)$, we have $q=B(\tilde{\Phi}_1, \lambda_1)$ with
\begin{align}   \label{Phi-6.1}
& \tilde{\Phi}_1 (x)    = (\frac{1}{\bar{\lambda}_1}-\frac{1}{\lambda_1})         \frac{1 }{ \gamma a^{(1)}(\lambda_1) }  e^{-i\lambda_1^2 x -2i \lambda_1^4 t} \varphi_{-}^{(1)} (x;\lambda_1) + (\frac{1}{\bar{\lambda}_1}-\frac{1}{\lambda_1}) \frac{1 }{ a^{(1)}(\lambda_1) }  e^{i\lambda_1^2 x -2i \lambda_1^4 t} \phi_{+}^{(1)} (x;\lambda_1),
\end{align}
where $a^{(1)} \neq 0$ by Proposition 4.3.

In addition, because $q^{(1)}(t, \cdot)$ exists for $t\in \mathbb{R}$, the corresponding Jost functions \\
 $\{ \varphi_{\pm}^{(1)}(t,x,\lambda), \phi_{\pm}^{(1)}(t,x,\lambda) \}$ exist for the time $t\in \mathbb{R}$.
 So we  denote $ \tilde{q}= B(\Phi_1, \lambda_1) q^{(1)}$ for  $t\in \mathbb{R}$. As $q^{(1)}(t,\cdot) = \tilde{q}(t,\cdot) \in Z_0$ for $t\in (-T,T)$ by uniqueness, the function $\tilde{q}$ is unique to the solution $q$ for the same Cauchy problem \eqref{DNLSII-eq-3}, which exists globally in time based on the bound estimate in Proposition 4.4.
It can be proved similarly that $ \| q^{(1)}(t,\cdot) \|_{H^2(\mathbb{R}) \cap H^{1,1}(\mathbb{R})} \leq M_T $ for time $t \in(-T,T)$, where $T>0$ is arbitrary and $M_T$ relies on $T$.
And combining Proposition 4.6, there exists $ \| q(t,\cdot) \|_{H^2(\mathbb{R}) \cap H^{1,1}(\mathbb{R})} \leq M_T $ for  $t \in(-T,T)$.
Therefore, the solution will not explode in a limited time and there exists a unique global solution $q(t,\cdot) \in Z_1$ for $t\in\mathbb{ R}$ with regard to the Cauchy problem \eqref{DNLSII-eq-3} for $q_0 \in Z_1\subset H^2(\mathbb{R}) \cap H^{1,1}(\mathbb{R})$.

By using the same discussion above and iterating $N$ times for the Darboux transformation, the global existence of $q(t,\cdot) \in Z_N \subset H^2(\mathbb{R}) \cap H^{1,1}(\mathbb{R})$  is proved from the global existence of $q^{(N)}(t,\cdot) \in Z_0 \subset H^2(\mathbb{R}) \cap H^{1,1}(\mathbb{R})$, $t\in \mathbb{R}$.  $\Box$

\section*{Acknowledgements}

This work was supported by the National Natural Science Foundation of China grant  No.12371256, No.11971475.\\

%
%
%
%



\end{CJK*}
\end{document}